\documentclass[12pt]{amsart}

\usepackage{amssymb}

\usepackage{enumitem}
\usepackage{hyperref}
\usepackage{cleveref}

\usepackage{graphicx}

\makeatletter
\@namedef{subjclassname@2020}{%
  \textup{2020} Mathematics Subject Classification}
\makeatother

\usepackage[T1]{fontenc}

\newtheorem*{theorem*}{Theorem}
\newtheorem{theorem}{Theorem}[section]
\newtheorem{corollary}[theorem]{Corollary}
\newtheorem{lemma}[theorem]{Lemma}
\newtheorem{proposition}[theorem]{Proposition}

\newtheorem{question}[theorem]{Question}

\theoremstyle{definition}

\newtheorem{remark}[theorem]{Remark}

\numberwithin{equation}{section}

\begin{document}


\baselineskip=17pt


\title[On cleanness of $AW^*$-algebras]{On cleanness of $AW^*$-algebras}

\author[L. Cui]{Lu Cui}
\email{nilu@upc.edu.cn}

\author[M. H. Ma]{Minghui Ma}
\email{minghuima@dlut.edu.cn}

\date{}

\begin{abstract}
  A ring is called clean if every element is the sum of an invertible element and an idempotent.
  This paper investigates the cleanness of $AW^*$-algebras.
  We prove that all finite $AW^*$-algebras are clean, affirmatively solving a question posed by Va\v{s}.
  We also prove that all countably decomposable infinite $AW^*$-factors are clean.
  A $*$-ring is called almost $*$-clean if every element can be expressed as the sum of a non-zero-divisor and a projection.
  We show that an $AW^*$-algebra is almost $*$-clean if and only if it is finite.
\end{abstract}



\maketitle

\section{Introduction}

Exchange rings are an important class of research objects in ring theory, which include (von Neumann) regular rings \cite{von} and have a close connection with $C^*$-algebras.
In order to characterize exchange rings in which all idempotents are in the center, Nicholson introduced the concept of clean rings and proved that clean rings are all exchange rings \cite{Nicholson}.
However, there exist exchange rings that are not clean rings \cite{Camillo3,Handelman}.

A ring is called \emph{clean (resp. almost clean)} if every element is the sum of an invertible element (resp. a non-zero-divisor) and an idempotent.
There are many examples of clean rings \cite{Camillo2,Nicholson}.
Clean rings are closely related to unit-regular rings.
A ring is a unit-regular ring if and only if every element can be written as the product of an idempotent and an invertible element.
Consequently, clean rings can be viewed as the additive analogues of unit-regular rings.
It is noteworthy that a unit-regular ring is inherently a clean ring \cite{Camillo1,Camillo3}.
However, a clean ring is not necessarily a unit-regular ring \cite{Nicholson}.
Furthermore, from the viewpoint of elements, a unit-regular element need not be clean \cite{Khurana}.

In 2005, Lam posed the question of which von Neumann algebras are clean as rings at the Conference on Algebra and Its Applications, held at Ohio University, Athens, OH.
Inspired by Lam's question, Va\v{s} studied the cleanness of Baer $*$-rings and von Neumann algebras \cite{Vas}.
Since every von Neumann algebra is a $*$-ring, Va\v{s} introduced a more natural definition of cleanness for $*$-rings.
A $*$-ring is called \emph{$*$-clean (resp. almost $*$-clean)} if every element can be expressed as the sum of an invertible element (resp. a non-zero-divisor) and a projection.
Furthermore, a ring (resp. $*$-ring) is referred to as \emph{strongly clean (resp. strongly $*$-clean)} if every elemnent is the sum of of an invertible element and an idempotent (resp. a projection) that commutes.

The exploration of cleanness in elements and rings has garnered significant attention.
We refer to the references \cite{Mary,Zhu} therein.
Recent advancements have extended the study of cleanness in von Neumann algebras.
In \cite{Cui}, the authors proved that all finite von Neumann algebras and separable infinite factors are clean.
Moreover, they showed that a von Neumann algebra $\mathcal{A}$ is strongly clean if and only if there exists a finite number of mutually orthogonal central projections $P_i$ with sum $I$ such that $P_i\mathcal{A}$ is of type $\mathrm{I}_{n_i}$ with $n_i<\infty$, and $\mathcal{A}$ is strongly $*$-clean if and only if $\mathcal{A}$ is abelian.

An $AW^*$-algebra is a generalization of a von Neumann algebra \cite{Kaplansky}.
An \emph{$AW^*$-algebra} is a $C^*$-algebra that satisfies the following conditions:
$(A)$ In the partially ordered set of projections, any set of orthogonal projections has a least upper bound ({\rm LUB});
$(B)$ Any maximal commutative self-adjoint subalgebra is generated by its projections.
It is well-known that von Neumann algebras are $AW^*$-algebras, but the converse is not true\cite{Maitland3}.
For instance, the algebra of all bounded Baire functions on the real line modulo the set of first category is an abelian $AW^*$-algebra, which is not a von Neumann algebra \cite{Dixmier,Saito}.
Despite this, $AW^*$-algebras share many nice properties of von Neumann algebras, such as spectral decomposition and polar decomposition \cite{Berberian,Frank,Maitland1}.
An $AW^*$-algebra of type I is a von Neumann algebra if and only if its center is a von Neumann algebra \cite{Kaplansky1}.
Consequently, any $AW^*$-factor of type I is a von Neumann algebra.
However, there exist $AW^*$-factors of type III that are not von Neumann algebras \cite{Maitland2,Maitland3}.

In \cite{Vas}, Va\v{s} investigated the cleanness of $AW^*$-algebras.
She demonstrated that a finite $AW^*$-algebra of type $\mathrm{I}$ is almost $*$-clean, a regular finite $AW^*$-algebra of type $\mathrm{I}$ is $*$-clean, and posed the following questions.
\begin{question}\label{question}
    Which finite $AW^*$-algebras of type $\mathrm{I}$ are clean?
    Are $AW^*$-algebras of type $\mathrm{II}_1$ $($almost$)$ clean?
\end{question}

\noindent
In 2013, Akalan and Va\v{s} studied \Cref{question} and proved that all finite $AW^*$-algebras are almost clean \cite{Akalan}.
In this paper, we strengthen their result and obtain the following theorems.

\begin{theorem*}
  {\rm\textbf{\ref{section 4}}}\quad
  All finite $AW^*$-algebras are clean.
  More precisely, let $T$ be an element in a finite $AW^*$-algebra $\mathcal{A}$, then there exists an idempotent $P\in\mathcal{A}$ such that $T-P$ is invertible and $\|(T-P)^{-1}\|\leqslant 4$.
\end{theorem*}

\begin{theorem*}
  {\rm\textbf{\ref{section 4-2}}}\quad
  An $AW^*$-algebra $\mathcal{A}$ is almost $*$-clean if and only if $\mathcal{A}$ is finite.
\end{theorem*}\noindent
Regarding the cleanness of infinite $AW^*$-algebras, we derive the following result.
\begin{theorem*}
  {\rm\textbf{\ref{section 5}}}\quad
    All countably decomposable infinite $AW^*$-factors are clean.
\end{theorem*}

This paper is structured as follows.
In \Cref{preliminaries}, we introduce some basics about $AW^*$-algebras, such as the left and right projections, spectral projections and Halmos' two projections theorem.
We present auxiliary results in \Cref{Auxiliary results}, which are instrumental in subsequent proofs.
In \Cref{finite AW-algebras}, we consider the cleanness of finite $AW^*$-algebras and establish \Cref{section 4} and \Cref{section 4-2}.
In \Cref{infinite AW-factors}, we focus on countably decomposable infinite $AW^*$-factors and provide the proof of \Cref{section 5}.
 
\section{Notation and preliminaries}\label{preliminaries}

We introduce fundamental concepts of $AW^*$-algebras, including left and right projections, spectral projections and Halmos' two projections theorem.
For the sake of completeness, we also include proofs for some results that are likely well-known to experts.

Throughout the paper, let $\mathcal{A}$ be an $AW^*$-algebra.
A projection $E$ in $\mathcal{A}$ is considered \emph{countably decomposable} relative to $\mathcal{A}$ if any  orthogonal family of nonzero subprojections of $E$ in $\mathcal{A}$ is countable.
We define $\mathcal{A}$ to be \emph{countably decomposable} if the identity $I$ is countably decomposable relative to $\mathcal{A}$.
An \emph{$AW^*$-factor} is an $AW^*$-algebra whose center consists of scalar multiples of the identity $I$.
In a countably decomposable $AW^*$-factor, all infinite projections are equivalent.

For any $T\in\mathcal{A}$, its right annihilator $\{S\in\mathcal{A}: TS=0\}$ is of the form $(I-E)\mathcal{A}$ for some projection $E\in\mathcal{A}$.
We refer to $E$ as the \emph{right projection} of $T$, denoted by $R(T)$.
Analogously, there is a dual definition of the \emph{left projection} of $T$, denoted by $L(T)$.
It is evident that $R(T)=E$ if and only if $TE=T$ and $TS=0$ implies $ES=0$ for any $S\in\mathcal{A}$.
Moreover, we have $R(T)=L(T^*)$.
For any projections $E,F\in\mathcal{A}$, we use the notation $E^{\perp}=I-E$ and observe that
\begin{equation*}
  R(EF^\perp)=(E\vee F)-F,\quad L(EF^\perp)=E-(E\wedge F).
\end{equation*}
We say that $E$ is \emph{equivalent} to $F$ if there exists a partial isometry $V\in\mathcal{A}$ such that $V^*V=E$ and $VV^*=F$, denoted by $E\sim F$.
For any $T\in\mathcal{A}$, we have $R(T)\sim L(T)$.
Consequently, Kaplansky's formula holds:
\begin{equation*}
    (E\vee F)-F\sim E-(E\wedge F).
\end{equation*}
For further details, please refer to \cite{Kaplansky}.

The following result asserts that if an element $T\in\mathcal{A}$ is left (or right) invertible in $L(T)\mathcal{A}R(T)$, then it is invertible in $L(T)\mathcal{A}R(T)$.

\begin{proposition}\label{left invertible}
  Let $T$ be an element in an $AW^*$-algebra $\mathcal{A}$.
  Suppose there exists $S_0\in\mathcal{A}$ such that $S_0T=R(T)$ or $TS_0=L(T)$.
  Then there exists $S\in\mathcal{A}$ such that $ST=R(T)$ and $TS=L(T)$.
\end{proposition}

\begin{proof}
Without loss of generality, we can assume that $S_0T=R(T)$.
Since $R(T)=T^*S_0^*S_0T\leqslant\|S_0\|^2 T^*T$, $T^*T$ is invertible in $R(T)\mathcal{A}R(T)$.
Let $S_1$ be the inverse of $T^*T$ in $R(T)\mathcal{A}R(T)$ and $S=S_1 T^*$.
Then $SL(T)=S$ and $ST=R(T)$.
  Since $TST=TR(T)=T$, we have $TSL(T)=L(T)$.
  Therefore, $TS=L(T)$.
\end{proof}

From the proof, the operator $S$ in \Cref{left invertible} actually lies in the $C^*$-algebra generated by $T$.

\begin{corollary}\label{zero-divisor}
  Let $V$ be a non-unitary isometry in an $AW^*$-algebra $\mathcal{A}$, i.e., $V^*V=I$ and $VV^*<I$.
  Then $L(V-A)<I$ for any $A\in\mathcal{A}$ with $\|A\|<1$.
\end{corollary}

\begin{proof}
    Since $\|A\|<1$, we obtain that $V^*(V-A)=I-V^*A$ is invertible.
    Let $S=(I-V^*A)^{-1}V^*$.
    Then we have $S(V-A)=I$, and hence $R(V-A)=I$.
    Suppose on the contrary that $L(V-A)=I$.
    Then $V-A$ is invertible by \Cref{left invertible}.
    Thus $V^*=(I-V^*A)(V-A)^{-1}$ is invertible.
    That is a contradiction.
\end{proof}

The following result can be readily obtained by examining the universal representation of $\mathcal{A}$.
However, we provide an intrinsic proof here.

\begin{lemma}\label{A-B}
Let $A,B\in\mathcal{A}$ and $a>b>0$.
If $E$ is a projection such that $\|BE\|\leqslant b$ and $EA^*AE\geqslant a^2E$, then $E(A-B)^*(A-B)E\geqslant(a-b)^2E$.
\end{lemma}

\begin{proof}
Since $|AE|^2=EA^*AE\geqslant a^2E$, $|AE|$ is invertible in $E\mathcal{A}E$, whose inverse is denoted by $S$.
Then $\|S\|\leqslant\frac{1}{a}$.
Let $AE=U|AE|$ be the polar decomposition, and let $S_1=SU^*$.
It is clear that $S_1AE=E$.
Since $\|S_1BE\|\leqslant\frac{b}{a}<1$, $E-S_1BE$ is invertible in $E\mathcal{A}E$, whose inverse is denoted by $S_2$.
Then $\|S_2\|\leqslant\frac{a}{a-b}$.
It follows that
\begin{equation*}
\begin{split}
    E & =(E-S_1BE)^*S_2^*S_2(E-S_1BE) \leqslant\|S_2\|^2(E-S_1BE)^*(E-S_1BE)\\
    & =\|S_2\|^2E(A-B)^*S_1^*S_1(A-B)E\leqslant\|S_2\|^2\|S_1\|^2E(A-B)^*(A-B)E.
\end{split}
\end{equation*}
This completes the proof.
\end{proof}

\begin{lemma}\label{TES}
  Let $E$ be a projection in $\mathcal{A}$ and $T\in\mathcal{A}$.
  If $ET^*TE\geqslant a^2E$ for some $a>0$, then there exists positive $S\in E\mathcal{A}E$ such that $TS$ is a partial isometry with initial projection $E$ and finial projection $L(TE)$, i.e., $(TS)^*TS=E$ and $TS(TS)^*=L(TE)$.
\end{lemma}

\begin{proof}
Let $A=TE$.
Since $A^*A\geqslant a^2E$, there exists positive $S\in E\mathcal{A}E$ such that
\begin{equation*}
    S^2A^*A=A^*AS^2=E.
\end{equation*}
Then we have $(AS)^*AS=E$ and $TS=AS$ is a partial isometry.
Note that $L(A)AS(AS)^*=AS(AS)^*$ and $AS(AS)^*A=A$.
Then we obtain $AS(AS)^*\leqslant L(A)\leqslant AS(AS)^*$.
\end{proof}

The spectral decomposition is a fundamental technique in the study of operator algebras.
Previous research has explored the spectral decomposition in $AW^*$-algebras \cite{Frank,Maitland1}.
In this paper, we only need to use the spectral projection of a self-adjoint element $T$ in $\mathcal{A}$ with respect to an interval.
This projection can be defined using the left projection of $T$.
Specifically, the \emph{spectral projection} of $T$ corresponding to the interval $(-\infty,c]$ is given by $P=L(f(T))^\perp$, where $f(x)=\max\{x-c,0\}$.
It is evident from the following lemma that if $TS=ST$, then $PS=SP$ for every $S\in\mathcal{A}$, implying in particular that $PT=TP$.
Additionally, we have the following properties:
\begin{equation*}
    PTP\leqslant cP,\quad P^\perp TP^\perp\geqslant cP^\perp.
\end{equation*}

\begin{lemma}
Let $c$ be a real number and $f(x)=\max\{x-c,0\}$ a real-valued continuous function on $\mathbb{R}$.
Let $T$ be a self-adjoint element in an $AW^*$-algebra $\mathcal{A}$ and $Q=L(f(T))=R(f(T))$.
Then $TQ\geqslant cQ$ and $TQ^{\perp}\leqslant cQ^{\perp}$.
Moreover, $TS=ST$ implies $QS=SQ$ for every $S\in\mathcal{A}$.
\end{lemma}

\begin{proof}
  Without loss of generality, we assume $c=0$.
  Since $(x-f(x))f(x)=0$, we have $(T-f(T))f(T)=0$ and hence $(T-f(T))Q=0$.
  It follows that $TQ=f(T)\geqslant 0$ and $TQ^\perp=T-f(T)\leqslant 0$.
  If $TS=ST$ for some $S\in\mathcal{A}$, then $f(T)S=Sf(T)$.
  Since $QSf(T)=Qf(T)S=f(T)S=Sf(T)$, we have $QSQ=SQ$.
  Similarly, since $f(T)SQ=Sf(T)Q=Sf(T)=f(T)S$, we have $QSQ=QS$.
Therefore, $QS=SQ$.
\end{proof}

Let $P$ and $Q$ be projections in $\mathcal{A}$.
We say that $P$ and $Q$ are in \emph{generic position} if
\begin{equation*}
    P\wedge Q=P\wedge Q^\perp
    =P^\perp\wedge Q
    =P^\perp\wedge Q^\perp=0.
\end{equation*}
If $P$ and $Q$ are projections in generic position, then by \cite[Lemma 5.3]{Kaplansky},
  \begin{equation*}
    R(QP)=Q\vee P^\perp-P^\perp=P,\quad L(QP)=Q-Q\wedge P^\perp=Q.
  \end{equation*}
Using a proof similar to Halmos' two projections theorem \cite{Halmos}, we can obtain the following result for $AW^*$-algebras.

\begin{proposition}\label{Halmos two projections}
  Let $P$ and $Q$ be projections in generic position in an $AW^*$-algebra $\mathcal{A}$.
  Then there exists a system of matrix units $\{E_{ij}\}_{1\leqslant i,j\leqslant 2}$ and an operator $H$ in $\mathcal{A}$, such that
  \begin{enumerate}
    \item [$(1)$] $E_{11}+E_{22}=I$;
    \item [$(2)$] $0\leqslant H\leqslant I$ and $L(H)=R(H)=I$;
    \item [$(3)$] $H$ commutes with $\{E_{ij}\}_{1\leqslant i,j\leqslant 2}$;
    \item [$(4)$] $P=E_{11}$ and $Q=E_{11}H+(E_{12}+E_{21})\sqrt{H(I-H)}+E_{22}(I-H)$.
  \end{enumerate}
In this case, we have $(P-Q)^2=I-H$.
\end{proposition}

\begin{remark}\label{P vee Q}
Using the notation in {\rm\Cref{Halmos two projections}}, it is easy to check that
  \begin{equation*}
    \begin{split}
      Q^\perp(I+iE_{21}-iE_{12})Q^\perp&=Q^\perp,\\
      \left[(I-iE_{21}+iE_{12})-2Q\right]^2&=\left[(I-iE_{21}+iE_{12})-2P\right]^2=2I.
    \end{split}
  \end{equation*}
\end{remark}

\begin{remark}\label{EFPQ}
  Let $E$, $F$ be projections in $\mathcal{A}$ and
  \begin{equation*}
    P=E-E\wedge F-E\wedge F^\perp,\quad Q=F-E\wedge F-E^\perp\wedge F,\quad I_0=P\vee Q.
  \end{equation*}
  Then $E\vee F=I_0+E\wedge F+E\wedge F^\perp+E^\perp\wedge F$.
  Moreover, $P$ and $Q$ are in generic position in $I_0\mathcal{A}I_0$.
  By {\rm\Cref{Halmos two projections}}, we can write
  \begin{equation*}
      P=E_{11},\quad Q=E_{11}H+(E_{12}+E_{21})\sqrt{H(I_0-H)}+E_{22}(I_0-H).
  \end{equation*}
  Suppose $E\wedge F^\perp\sim E^\perp\wedge F$, i.e., $U^*U=E\wedge F^\perp$ and $UU^*=E^\perp\wedge F$ for some partial isometry $U\in\mathcal{A}$.
  Let
  \begin{equation*}
    P_0=\frac{1}{2}(I_0+iE_{21}-iE_{12})+E\wedge F+\frac{1}{2}(E\wedge F^{\perp}+E^{\perp}\wedge F+U+U^*).
  \end{equation*}
  Then $P_0$ is a projection.
  By {\rm\Cref{P vee Q}}, $F^\perp P_0 F^\perp=\frac{1}{2}(I_0-Q)+\frac{1}{2}E\wedge F^{\perp}$,
  \begin{equation*}
     (P_0-E^\perp)^2=(P_0-F^\perp)^2=\frac{1}{2}I_0+E\wedge F+\frac{1}{2}(E\wedge F^\perp+E^\perp\wedge F)+E^\perp\wedge F^\perp.
  \end{equation*}
  It follows that $\|P_0 F^\perp\|=\frac{1}{\sqrt{2}}$, $P_0-E^\perp$ and $P_0-F^\perp$ are invertible in $\mathcal{A}$.
\end{remark}

\section{Auxiliary results}\label{Auxiliary results}

This section presents some auxiliary results that will be utilized in subsequent proofs.
Let $\mathcal{F}_\mathcal{A}$ be the two-sided ideal generated by finite projections in $\mathcal{A}$, and $\mathcal{K}_\mathcal{A}$ the norm-closure of $\mathcal{F}_\mathcal{A}$ in $\mathcal{A}$.
Operators in $\mathcal{F}_\mathcal{A}$ and $\mathcal{K}_\mathcal{A}$ are called \emph{finite and compact relative to $\mathcal{A}$}, respectively.
It is clear that $T\in\mathcal{F}_\mathcal{A}$ if and only if there exists a finite projection $E\in\mathcal{A}$ such that $T\in E\mathcal{A}E$.
The following result elucidates certain properties of operators in $\mathcal{K}_{\mathcal{A}}$.

\begin{proposition}\label{compact equivalent}
  Let $T$ be an element in an infinite $AW^*$-factor $\mathcal{A}$.
  Then the following statements are equivalent:
  \begin{enumerate}
    \item [$(1)$] $T$ is compact relative to $\mathcal{A}$.
    \item [$(2)$] There is no infinite projection $Q\in\mathcal{A}$ such that $TS=Q$ for some $S\in\mathcal{A}$.
    \item [$(3)$] The spectral projection of $|T|$ with respect to $(c,+\infty)$ is finite for any $c\in(0,1)$.
  \end{enumerate}
\end{proposition}

\begin{proof}
$(1)\Rightarrow(2)$:
Let $Q$ be a projection in $\mathcal{A}$ such that $TS=Q$ for some $S\in\mathcal{A}$.
Since $T$ is compact, we have $\|A-T\|<\frac{1}{1+\|S\|}$ for some $A\in\mathcal{F}_\mathcal{A}$.
It follows that $\|QASQ-Q\|=\|Q(A-T)SQ\|<1$.
Then $QASQ$ is invertible in $Q\mathcal{A}Q$.
Hence $Q=L(QASQ)\leqslant L(QA)\sim R(QA)\leqslant R(A)$.
Therefore, $Q$ is finite.

$(2)\Rightarrow(3)$: Let $T=V|T|$ be the polar decomposition of $T$ and
$P$ the spectral projection of $|T|$ with respect to $(c,+\infty)$.
Since $|T|P\geqslant cP$, there exists $S\in P\mathcal{A}P$ such that $|T|PS=P$.
Then $TPSV^*=VPV^*$.
Since $P\leqslant R(|T|)=V^*V$, we have $VPV^*\sim P$.
By condition $(2)$, $VPV^*$ is a finite projection.
Hence $P$ is finite.

  $(3)\Rightarrow(1)$: Let $P_n$ be the spectral projection of $|T|$ with respect to $[0,\frac{1}{n}]$ for any positive integer $n$.
  Then $P_n^\perp$ is finite and $P_nT^*TP_n\leqslant\frac{1}{n^2}P_n$.
  Hence $\|T-TP_n^\perp\|\leqslant\frac{1}{n}$.
  Since $TP_n^\perp\in\mathcal{F}_\mathcal{A}$, we obtain $T\in\mathcal{K}_\mathcal{A}$.
\end{proof}

\begin{corollary}\label{ET compact}
  Let $E$ be an infinite projection in an infinite $AW^*$-factor $\mathcal{A}$ and $T\in\mathcal{A}$.
  If $ET\in\mathcal{K}_\mathcal{A}$, then there is no finite projection $F\in\mathcal{A}$ such that $TS=F^\perp$ for some $S\in\mathcal{A}$.
\end{corollary}

\begin{proof}
Let $F$ be a projection such that $TS=F^\perp$ for some $S\in\mathcal{A}$.
We have $E\wedge F^\perp=(E\wedge F^\perp)ETS$.
Since $ET\in\mathcal{K}_{\mathcal{A}}$, $E\wedge F^\perp$ is finite by \Cref{compact equivalent}.
Since $E-E\wedge F^\perp\sim E\vee F^\perp-F^\perp\leqslant F$, $F$ is infinite.
This completes the proof.
\end{proof}

Let $\pi\colon\mathcal{A}\to\mathcal{A}/\mathcal{K}_{\mathcal{A}}$ be the quotient map.
The following result is an immediate consequence of \Cref{compact equivalent}.

\begin{lemma}\label{center}
Let $\mathcal{A}$ be an infinite $AW^*$-factor.
Then the center of $\mathcal{A}/\mathcal{K}_{\mathcal{A}}$ is trivial.
\end{lemma}

\begin{proof}
Let $\mathcal{Z}$ be the center of $\mathcal{A}/\mathcal{K}_{\mathcal{A}}$.
Then $\mathcal{Z}$ is a commutative $C^*$-algebra.
Suppose on the contrary that $\mathcal{Z}$ is non-trivial.
By continuous function calculus, there are nonzero elements $\pi(T)$ and $\pi(S)$ in $\mathcal{Z}$ with $\pi(T)\pi(S)=0$.
By \Cref{compact equivalent}, there exist $T_1,S_1\in\mathcal{A}$ such that $TT_1=P$ and $S_1S=Q$ are infinite projections in $\mathcal{A}$.
By \cite[Lemma 3.4]{Kaplansky}, we can assume $P\preceq Q$.
Then there is a partial isometry $V$ such that $VV^*=P$ and $V^*V\leqslant Q$.
Let $A=T_1VS_1$.
Then $TAS=V\notin\mathcal{K}_{\mathcal{A}}$.
Since $\pi(T),\pi(S)\in\mathcal{Z}$, we have $\pi(TAS)=\pi(T)\pi(A)\pi(S)=\pi(T)\pi(S)\pi(A)=0$.
That is a contradiction.
\end{proof}

In the following proposition, we characterize the operators that can be expressed as the sum of a scalar and a compact operator in an infinite $AW^*$-factor.

\begin{proposition}\label{T-cI compact}
  Let $T$ be an element in an infinite $AW^*$-factor $\mathcal{A}$.
  Then the following statements are equivalent:
  \begin{enumerate}
    \item [$(1)$] There exists $z\in\mathbb{C}$ such that $T-zI\in\mathcal{K}_\mathcal{A}$.
    \item [$(2)$] For any unitary operator $W\in\mathcal{A}$, $T-W^*TW\in\mathcal{K}_\mathcal{A}$.
    \item [$(3)$] For any projection $E\in\mathcal{A}$ with $E\sim E^\perp$, $E^\perp TE\in\mathcal{K}_\mathcal{A}$.
  \end{enumerate}
\end{proposition}

\begin{proof}
  $(1)\Rightarrow(2)$ and $(1)\Rightarrow(3)$ are clear.

  $(2)\Rightarrow(1)$: Since $WT-TW\in\mathcal{K}_{\mathcal{A}}$ for any unitary operator $W\in\mathcal{A}$, we have $AT-TA\in\mathcal{K}_{\mathcal{A}}$ for any $A\in\mathcal{A}$.
  Thus, $\pi(T)$ lies in the center of $\mathcal{A}/\mathcal{K}_{\mathcal{A}}$.
  By \Cref{center}, $T-zI\in\mathcal{K}_{\mathcal{A}}$ for some $z\in\mathbb{C}$.

  $(3)\Rightarrow(1)$: Let $E$ be a projection in $\mathcal{A}$ with $E\preceq E^{\perp}$.
  Since $E^{\perp}$ is infinite, by \cite[Lemma 4.5]{Kaplansky}, we can write $E^{\perp}=E_1+E_2$, where $E_1\sim E_2\sim E^{\perp}$.
  For $j=1,2$, let $F_j=E+E_j$.
  Then $F_j\sim F_j^{\perp}$ and $F_j^{\perp}TF_j, F_jTF_j^{\perp}\in\mathcal{K}_{\mathcal{A}}$.
  It follows that $\pi(T)$ commutes with $\pi(F_j)$.
  Therefore, $\pi(T)$ commutes with $\pi(E)=\pi(F_1)\pi(F_2)$.
  Thus, $\pi(T)$ lies in the center of $\mathcal{A}/\mathcal{K}_{\mathcal{A}}$.
  By \Cref{center}, $T-zI\in\mathcal{K}_{\mathcal{A}}$ for some $z\in\mathbb{C}$.
\end{proof}

The proof of the following result employs the fact that $F\wedge E^\perp=0$ implies $F\preceq E$ for projections $E,F\in\mathcal{A}$, which can be derived from Kaplansky's formula.

\begin{proposition}\label{|T*|-|VT*V-1|}
  Let $T$ be an element in an $AW^*$-algebra $\mathcal{A}$.
  Then the following statements are equivalent:
  \begin{enumerate}
    \item [$(1)$] There is a finite projection $Q\in\mathcal{A}$ such that $TS=Q^\perp$ for some $S\in\mathcal{A}$.
    \item [$(2)$] The spectral projection of $|T^*|$ with respect to $[0,c]$ is finite for some $c\in(0,1)$.
    \item [$(3)$] For any invertible element $V\in\mathcal{A}$, the spectral projection of $|(V^{-1}TV)^*|$ with respect to $[0,c_V]$ is finite for some $c_V\in(0,1)$.
  \end{enumerate}
\end{proposition}

\begin{proof}
  $(1)\Rightarrow(2)$: Since $TS=Q^\perp$, we have $Q^\perp=TSS^*T^*\leqslant\|S\|^2 TT^*$.
  Let $c=\min\left\{\frac{1}{2\|S\|},\frac{1}{2}\right\}$ and $P$ the spectral projection of $|T^*|$ with respect to $[0,c]$.
  Then $4c^2 P\wedge Q^\perp\leqslant(P\wedge Q^\perp)TT^*(P\wedge Q^\perp)\leqslant c^2 P\wedge Q^\perp$.
  It follows that $P\wedge Q^\perp=0$ and hence $P\preceq Q$.
  Therefore, $P$ is a finite projection.

  $(1)\Rightarrow(3)$: For any invertible element $V\in\mathcal{A}$, let
  \begin{equation*}
    T_0=V^{-1}TV,\quad S_0=V^{-1}SV,\quad  Q_0=V^{-1}QV,\quad P_0=V^{-1}Q^{\perp}V.
  \end{equation*}
  Since $T_0S_0P_0=P_0$, we have $T_0S_0L(P_0)=L(P_0)$.
  By the fact $(1)\Rightarrow(2)$, we need to show that $L(P_0)^{\perp}$ is finite.
  Since $P_0+Q_0=I$ and $Q_0R(Q_0)^{\perp}=0$, we have $P_0R(Q_0)^{\perp}=R(Q_0)^{\perp}$.
  Then $R(Q_0)^{\perp}=L(P_0R(Q_0)^{\perp})\leqslant L(P_0)$.
  Therefore, $L(P_0)^{\perp}\leqslant R(Q_0)\leqslant R(QV)\sim L(QV)\leqslant Q$ and $L(P_0)^{\perp}$ is finite.

   $(3)\Rightarrow(2)$: It is clear.

   $(2)\Rightarrow(1)$: Let $P$ be the spectral projection of $|T^*|$ with respect to $[0,c]$, then $TT^*P^\perp\geqslant c^2 P^\perp$ and hence $TT^*P^\perp$ is invertible in $P^\perp\mathcal{A}P^\perp$.
  Therefore, there exists $S_0\in P^\perp\mathcal{A}P^\perp$ such that $TT^*S_0=P^\perp$.
  Now we set $Q=P$ and $S=T^*S_0\in\mathcal{A}$.
  Then we complete the proof.
\end{proof}

\begin{lemma}\label{E-F invertible}
  Suppose $E$ and $F$ are two projections in an $AW^*$-algebra $\mathcal{A}$.
  Then the following statements are equivalent:
  \begin{itemize}
    \item[$(1)$] $E-F$ is invertible in $(E\vee F)\mathcal{A}(E\vee F)$;
    \item[$(2)$] $E=S(E\vee F-F)E$ and $E\vee F-F=(E\vee F-F)ES$ for some $S\in E\mathcal{A}(E\vee F-F)$;
    \item[$(3)$] $E\wedge F=0$ and $E\vee F=EA+FB$ for some $A,B\in\mathcal{A}$;
    \item[$(4)$] $E\wedge F=0$ and $E\vee F-F=(E\vee F-F)EA$ for some $A\in\mathcal{A}$.
  \end{itemize}
  If the conditions are satisfied, then $\|(E-F)|_{E\vee F}^{-1}\|=(1-\|EF\|^2)^{-\frac{1}{2}}=\|S\|$.
\end{lemma}

\begin{proof}
  We can assume $E\vee F=I$.
  It is clear that $(2)\Rightarrow(3)\Rightarrow(4)$.

  $(4)\Rightarrow(2)$: Since $L(F^\perp E)=F^\perp$ and $R(F^\perp E)=F^\perp\vee E^\perp-E^\perp=E$, the conclusion holds by \Cref{left invertible}.

  $(1)\Rightarrow(2)$: Note that $F^\perp E=(E-F)E=F^\perp(E-F)$.
   Then we have $SF^\perp E=E$ and $F^\perp ES=F^\perp$ for $S=E(E-F)^{-1}F^\perp\in E\mathcal{A}F^\perp$.

  $(2)\Rightarrow(1)$:
  It is clear that $E\wedge F=0$.
  Let $P=E-E\wedge F^\perp$, $Q=F-E^\perp\wedge F$ and $I_0=P\vee Q$.
  Then $P$ and $Q$ are projections in $I_0\mathcal{A}I_0$ in generic position, and $I=I_0+E\wedge F^\perp+E^\perp\wedge F$.
  By \Cref{Halmos two projections}, we can write
  \begin{equation*}
      P=E_{11},\quad Q=E_{11}H+(E_{12}+E_{21})\sqrt{H(I_0-H)}+E_{22}(I_0-H).
  \end{equation*}
  Since $E_{11}=SF^\perp EE_{11}=S(I_0-Q)E_{11}$, we have
  \begin{equation*}
      E_{11}=E_{11}(I_0-Q)S^*S(I_0-Q)E_{11}\leqslant\|S\|^2E_{11}(I_0-Q)E_{11}=\|S\|^2E_{11}(I_0-H).
  \end{equation*}
  It follows that $\|S\|^2(I_0-H)\geqslant I_0$ and hence $I_0-H$ is invertible in $I_0\mathcal{A}I_0$, whose inverse is denoted by $(I_0-H)^{-1}$.
  Then $P-Q$ is invertible in $I_0\mathcal{A}I_0$ with $(P-Q)^{-1}=E_{11}-(E_{12}+E_{21})\sqrt{H(I_0-H)^{-1}}-E_{22}$.
  Therefore, $E-F$ is invertible with $(E-F)^{-1}=(P-Q)^{-1}+E\wedge F^\perp-E^\perp\wedge F$.
  We complete the proof of the equivalence of (1), (2), (3) and (4).

  Under these conditions, $S=E_{11}-E_{12}\sqrt{H(I_0-H)^{-1}}+E\wedge F^\perp$ and $SS^*=E_{11}(I_0-H)^{-1}+E\wedge F^\perp$.
  Since $(E-F)^{-2}=(I_0-H)^{-1}+E\wedge F^\perp+E^\perp\wedge F$ and $EFE=E_{11}H$, we have
  \begin{equation*}
    \|(E-F)^{-1}\|=\|S\|=\sqrt{\|(I_0-H)^{-1}\|}=(1-\|H\|)^{-\frac{1}{2}}=(1-\|EF\|^2)^{-\frac{1}{2}}.
  \end{equation*}
  We complete the proof.
\end{proof}

\begin{corollary}\label{T invertible in A}
  Let $T$ be an element in an $AW^*$-algebra $\mathcal{A}$.
  Then there exists $S\in\mathcal{A}$ such that $TS=L(T)$ and $ST=I$ if and only if there is a projection $E\in\mathcal{A}$ satisfying the following conditions:
  \begin{enumerate}
    \item [$(a)$] $ET^*TE\geqslant a_1^2E$ and $E^\perp T^*TE^\perp\geqslant a_2^2E^\perp$ for some $a_1,a_2>0$;
    \item [$(b)$] $L(TE)-L(TE^\perp)$ is invertible in $L(T)\mathcal{A}L(T)$.
  \end{enumerate}
  If the conditions $(a)$ and $(b)$ are satisfied, we have
\begin{equation}\label{norm}
    \frac{1}{\|T\|}\leqslant\|S\|\sqrt{1-\|L(TE)L(TE^\perp)\|}\leqslant\frac{1}{\min\{a_1, a_2\}}.
\end{equation}
\end{corollary}

\begin{proof}
  Necessity:
  We have $E=ET^*S^*STE\leqslant\|S\|^2 ET^*TE$ for any projection $E\in\mathcal{A}$.
  Similarly, $E^\perp\leqslant\|S\|^2 E^\perp T^*TE^\perp$.
  Hence condition $(a)$ holds for $a_1=a_2=\|S\|^{-1}$.
  By \Cref{TES}, we have $L(TE)=TES_1$ and $L(TE^\perp)=TE^\perp S_2$ for some $S_1,S_2\in\mathcal{A}$.
  Let $P=L(TE)\wedge L(TE^\perp)$.
  Then $P=TES_1P=TE^\perp S_2P$.
  Since $SP=ES_1P=E^\perp S_2P$, we have $SP=0$ and $P=L(T)P=TSP=0$.
  Clearly, $L(T)=L(TE)\vee L(TE^\perp)$.
  Let $A=TES$ and $B=TE^\perp S$.
  Then $L(T)=L(TE)A+L(TE^\perp)B$.
  Hence $(b)$ holds by \Cref{E-F invertible}.

  Sufficiency: Suppose $TP=0$ for some projection $P\in\mathcal{A}$.
  Then we have $TEP=-TE^\perp P$.
  Since $L(TEP)=L(TE^\perp P)\leqslant L(TE)\wedge L(TE^\perp)=0$,
  we have $TEP=TE^\perp P=0$.
  By condition $(a)$, we have $EP=E^\perp P=0$, i.e., $P=0$.
  It follows that $R(T)=I$.
  By \Cref{TES}, we have $L(TE)=TES_1$ and $L(TE^\perp)=TE^\perp S_2$ for some $S_1,S_2\in\mathcal{A}$.
  Moreover, by \Cref{E-F invertible}, $L(T)=L(TE)A+L(TE^\perp)B$ for some $A,B\in\mathcal{A}$.
  Let $S_0=ES_1A+E^\perp S_2B$.
  Then we have $TS_0=L(T)$.
  By \Cref{left invertible}, there exists $S\in\mathcal{A}$ such that
  $ST=R(T)=I$ and $TS=L(T)$.

  Next we prove \eqref{norm}.
  Let $T_1=TE$, $T_2=TE^\perp$, $P_j=L(T_j)$ for $j=1,2$ and $\lambda=\|P_1P_2\|$.
  Since
    \begin{equation*}
  \begin{split}
      \lambda^2(T_1^*T_1+T_2^*T_2)-\lambda(T_1^*T_2+T_2^*T_1) & \geqslant T_1^*P_2T_1+\lambda^2T_2^*T_2-\lambda(T_1^*T_2+T_2^*T_1)\\ & =(T_1-\lambda T_2)^*P_2(T_1-\lambda T_2)\geqslant 0,
    \end{split}
  \end{equation*}
  we have $-\lambda(T_1^*T_1+T_2^*T_2)\leqslant T_1^*T_2+T_2^*T_1\leqslant \lambda(T_1^*T_1+T_2^*T_2)$.
  Therefore,
  \begin{align*}
      T^*T=(T_1^*T_1+T_2^*T_2)+(T_1^*T_2+T_2^*T_1)&\geqslant(1-\lambda)(T_1^*T_1+T_2^*T_2)\\
      &\geqslant(1-\lambda)\min\{a_1^2,a_2^2\}I.
  \end{align*}
  Then we obtain $(1-\lambda)\min\{a_1^2,a_2^2\}S^*S\leqslant I$. 
  Hence
  $\|S\|\sqrt{1-\lambda}\leqslant\frac{1}{\min\{a_1,a_2\}}$.

  Let $Q_1=P_1-P_1\wedge P_2^\perp$, $Q_2=P_2-P_1^\perp\wedge P_2$ and $I_0=Q_1\vee Q_2$.
  It follows that $Q_1$ and $Q_2$ are projections in generic position in $I_0\mathcal{A}I_0$, and $L(T)=I_0+P_1\wedge P_2^\perp+P_1^\perp\wedge P_2$.
  By \Cref{Halmos two projections}, we can write
  \begin{equation*}
      Q_1=E_{11},\quad Q_2=E_{11}H+(E_{12}+E_{21})\sqrt{H(I_0-H)}+E_{22}(I_0-H).
  \end{equation*}
  Let $V=E_{11}\sqrt{H}+E_{21}\sqrt{I_0-H}$.
  Then $V^*V=Q_1$ and $VV^*=Q_2$.
  By \Cref{TES}, we have $P_1=TES_1$ and $P_2=TE^\perp S_2$ for some $S_1,S_2\in\mathcal{A}$.
  Then $E_{11}=TA_1$ and $V=TA_2$, where $A_1=ES_1E_{11}$ and $A_2=E^\perp S_2V$.
  Hence
  \begin{equation*}
  \begin{split}
      A_1^*A_1+A_2^*A_2 & =(A_1-A_2)^*(A_1-A_2)=(A_1-A_2)^*T^*S^*ST(A_1-A_2)\\
      & \leqslant\|S\|^2(A_1-A_2)^*T^*T(A_1-A_2)=2\|S\|^2E_{11}(I_0-\sqrt{H}).
  \end{split}
  \end{equation*}
  Note that $\lambda=\|\sqrt{H}\|$.
  Then for any $\varepsilon>0$, there exist nonzero projections $F_2\leqslant F_1\leqslant F_0\leqslant E_{11}$ such that
  \begin{equation*}
  \begin{split}
      F_0\sqrt{H}F_0 & \geqslant(\lambda-\varepsilon)F_0,\\
      F_1A_1^*A_1F_1 & =F_1(F_0A_1^*A_1F_0)F_1\geqslant(\|A_1F_0\|^2-\varepsilon)F_1,\\
      F_2A_2^*A_2F_2 & \geqslant F_2(F_1A_2^*A_2F_1)F_2\geqslant(\|A_2F_1\|^2-\varepsilon)F_2.
  \end{split}
  \end{equation*}
  It follows that
  \begin{equation*}
  \begin{split}
      2\|S\|^2(1-\lambda+\varepsilon)F_2 & \geqslant 2\|S\|^2F_2(I_0-\sqrt{H})F_2\geqslant F_2(A_1^*A_1+A_2^*A_2)F_2\\ & \geqslant(\|A_1F_0\|^2+\|A_2F_1\|^2-2\varepsilon)F_2.
  \end{split}
  \end{equation*}
  Since $F_0=TA_1F_0$ and $F_1=V^*TA_2F_1$, we have $\|T\|\|A_1F_0\|\geqslant 1$ and $\|T\|\|A_2F_1\|\geqslant 1$.
  It follows that $\|S\|^2(1-\lambda+\varepsilon)\geqslant\frac{1}{\|T\|^2}-\varepsilon$, and hence $\|S\|\sqrt{1-\lambda}\geqslant\frac{1}{\|T\|}$.
\end{proof}

\begin{lemma}\label{idempotent}
  Let $T$ be an element in an $AW^*$-algebra $\mathcal{A}$.
  Suppose $E$ is a projection in $\mathcal{A}$ such that $E-ETE$ is invertible in $E\mathcal{A}E$.
  Then for every $A\in E^\perp\mathcal{A}E$, there exists an idempotent $P\in\mathcal{A}$ such that
  \begin{enumerate}
    \item [$(1)$] $R(P)=E$,
    \item [$(2)$] $P-E\in E^\perp\mathcal{A}E$,
    \item [$(3)$] $\|P\|\leqslant (1+\|A\|)(1+\|TE\|)$,
    \item [$(4)$] $L((T-P)E)=L(E+A)$.
  \end{enumerate}
\end{lemma}

\begin{proof}
  Let $P=E+E^\perp TE+A(E-ETE)\in\mathcal{A}$.
  Then $P$ is an idempotent satisfying $(1)$, $(2)$ and $(3)$.
  Note that $P=TE+(E+A)(E-ETE)$ and $E-ETE$ is invertible in $E\mathcal{A}E$.
  We get $(T-P)E=-(E+A)(E-ETE)$ and $L((T-P)E)=L(E+A)$.
\end{proof}

\section{Cleanness of finite $AW^*$-algebras}\label{finite AW-algebras}

In this section, we utilize the results from the previous section to establish the cleanness of finite $AW^*$-algebras.

\begin{lemma}\label{Tclean}
  Let $T$ be an element in an $AW^*$-algebra $\mathcal{A}$.
  Suppose there exists a projection $E\in\mathcal{A}$ such that
  \begin{itemize}
    \item[$(1)$] $E-ETE$ is invertible in $E\mathcal{A}E$ and $E^\perp T^*TE^\perp\geqslant c^2 E^\perp$ for some $c>0$,
    \item[$(2)$] $E\wedge L(TE^\perp)\sim E^\perp\wedge(I-L(TE^\perp))$.
  \end{itemize}
  Then there is an idempotent $P\in\mathcal{A}$ such that $\|P\|\leqslant 2+2\|TE\|$ and $T-P$ is invertible with
  \begin{equation*}
    \|(T-P)^{-1}\|\leqslant\frac{2}{\min\{\|(E-ETE)^{-1}\|^{-1}, c\}}.
  \end{equation*}
\end{lemma}

\begin{proof}
  Let $F=I-L(TE^\perp)$.
  Then $E\wedge F^\perp\sim E^\perp\wedge F$.
  We will use the notation introduced in \Cref{EFPQ}.
  Let $A=iE_{21}+U\in E^\perp\mathcal{A}E$ and $P$ the idempotent given by \Cref{idempotent}.
  Then we have
  \begin{equation*}
    L((T-P)E)=L(E+A)=L((E+A)(E+A)^*)=P_0.
  \end{equation*}
  Note that $L((T-P)E^\perp)=L(TE^\perp)=F^\perp$.
  By \Cref{EFPQ}, we obtain that $L((T-P)E)-L((T-P)E^\perp)$ is invertible in $\mathcal{A}$.
  Thus $T-P$ satisfies \Cref{T invertible in A} $(b)$.
  Since $E(T-P)E=-(E-ETE)$, we have
  \begin{equation*}
      E(T-P)^*(T-P)E
      \geqslant E(T-P)^*E(T-P)E=(E-ETE)^*(E-ETE).
  \end{equation*}
  Then $T-P$ satisfies \Cref{T invertible in A} $(a)$ for $a_1=\frac{1}{\|(E-ETE)^{-1}\|}$ and $a_2=c$.
  Hence, $T-P$ is invertible in $\mathcal{A}$ by \Cref{T invertible in A} and $\|L((T-P)E)F^\perp\|=\frac{1}{\sqrt{2}}$ by \Cref{EFPQ}.
  We complete the proof.
\end{proof}

\begin{theorem}\label{section 4}
  All finite $AW^*$-algebras are clean.
  More precisely, let $T$ be an element in a finite $AW^*$-algebra $\mathcal{A}$, then there exists an idempotent $P\in\mathcal{A}$ such that $T-P$ is invertible and $\|(T-P)^{-1}\|\leqslant 4$.
\end{theorem}

\begin{proof}
  If $\|T\|\leqslant\frac{1}{2}$, then $I-T$ is invertible and $\|(I-T)^{-1}\|\leqslant 2$.
  We only need to consider the case $\|T\|>\frac{1}{2}$.
  Let $E$ be the spectral projection of $|T|$ with respect to $[0,\frac{1}{2}]$.
  Then $E^\perp T^*TE^\perp\geqslant \frac{1}{4}E^\perp$ and  $\|ETE\|\leqslant\|TE\|\leqslant\frac{1}{2}$.
  Hence, $R(TE^\perp)=E^\perp$, $E-ETE$ is invertible in $E\mathcal{A}E$ and   $\|(E-ETE)^{-1}\|\leqslant 2$.
  Let $F=I-L(TE^\perp)$.
  Then $E^\perp=R(TE^\perp)\sim L(TE^\perp)=F^\perp$.
  By Kaplansky's formula, we have $E^\perp-E^\perp\wedge F\sim E^\perp\vee F-F=F^\perp-E\wedge F^\perp$.
  Since $\mathcal{A}$ is finite, $E\wedge F^\perp\sim E^\perp\wedge F$.
  We complete the proof by \Cref{Tclean}.
\end{proof}

\begin{theorem}\label{section 4-2}
  An $AW^*$-algebra $\mathcal{A}$ is almost $*$-clean if and only if $\mathcal{A}$ is finite.
\end{theorem}

\begin{proof}
  Necessity:
  Suppose $\mathcal{A}$ is not finite. Without loss of generality, we can assume $\mathcal{A}$ is properly infinite.
  By \cite[Lemma 4.5]{Kaplansky}, there is a family of orthogonal projections $\{E_n\}_{n=1}^\infty$ in $\mathcal{A}$ such that $E_n\sim I$ for each $n\geqslant 1$ and $\sum_{n=1}^{\infty}E_n=I$.
  Let $V_n$ be a partial isometry in $\mathcal{A}$ such that $V_n^*V_n=E_n$ and $V_n V_n^*=E_{n+1}$.
  By \cite[\S 20, Theorem 1]{Berberian}, we define $V=\sum_{n=1}^{\infty}V_n\in\mathcal{A}$.
  Then $V^*V=I$ and $VV^*=I-E_1<I$.
  By \Cref{zero-divisor}, we have $L(2V-P)<I$ for any projection $P\in\mathcal{A}$.
  Hence $2V$ is not almost $*$-clean.
  Therefore, $\mathcal{A}$ is not almost $*$-clean.

  Sufficiency:
  Suppose $\mathcal{A}$ is finite.
  Let $T$ be an operator in $\mathcal{A}$, $E=L(T)^\perp$ and $F=R(T)^\perp$.
  Then $E^\perp\sim F^\perp$.
  Using a similar proof of \Cref{section 4}, we can obtain $E\wedge F^\perp\sim E^\perp\wedge F$.
  We will use the notation introduced in \Cref{EFPQ}.
  It follows that $F-P_0^\perp$ and $P_0-E^\perp$ are invertible in $\mathcal{A}$.
  By \Cref{E-F invertible}, $L(T)\wedge P_0=R(T)^\perp\wedge P_0^\perp=0$.
  Let $P_1=R(T-P_0)$.
  Then we have $TP_1^\perp=P_0P_1^\perp$ and hence $L(TP_1^\perp)=L(P_0P_1^\perp)\leqslant L(T)\wedge P_0=0$.
  We obtain $TP_1^\perp=P_0P_1^\perp=0$.
  Then $P_1^\perp\leqslant R(T)^\perp\wedge P_0^\perp=0$.
  Therefore, $L(T-P_0)=R(T-P_0)=I$.
  Then $T-P_0$ is a non-zero-divisor.
  Thus, $T$ is almost $*$-clean.
  Therefore, $\mathcal{A}$ is almost $*$-clean.
\end{proof}

Naturally, we have the following problem.

\begin{question}
    Are all finite $AW^*$-algebras $*$-clean?
\end{question}

\section{Cleanness of countably decomposable infinite $AW^*$-factors}\label{infinite AW-factors}

In this section, we prove that all countably decomposable infinite $AW^*$-factors are clean.

\begin{lemma}\label{T and I-T}
  Let $T$ be an element in an $AW^*$-algebra $\mathcal{A}$.
  If the spectral projection of $|T|$ with respect to $(c,+\infty)$ is finite for some $c\in(0,1)$, then the spectral projection of $|I-T|$ with respect to $[0,b]$ is finite for any $b\in[0,1-c)$.
\end{lemma}

\begin{proof}
Let $E$ be the spectral projection of $|T|$ with respect to $[0,c]$ and $F$ the spectral projection of $|I-T|$ with respect to $[0,b]$.
Then $ET^*TE\leqslant c^2 E$ and $F(I-T)^*(I-T)F\leqslant b^2 F$.
It follows that
\begin{equation*}
    \|E\wedge F\|\leqslant\|T(E\wedge F)\|+\|(I-T)(E\wedge F)\|\leqslant\|TE\|+\|(I-T)F\|\leqslant c+b<1.
\end{equation*}
Therefore, we have $E\wedge F=0$ and $F=F-F\wedge E\sim F\vee E-E\leqslant E^\perp$.
Since $E^\perp$ is finite, $F$ is also finite.
\end{proof}

\begin{lemma}\label{TF<a}
  Let $E$ be a projection in an $AW^*$-algebra $\mathcal{A}$, $T\in\mathcal{A}$ and $F$ the spectral projection of $|T|$ with respect to $[0,c]$ for some $c>0$.
  If there exists $A\in\mathcal{A}$ with $\|A\|<c$ and $\|(T+A)E\|<c-\|A\|$, then $E\preceq F$.
  In particular, if $\|TE\|<c$, then $E\preceq F$.
\end{lemma}

\begin{proof}
  Since $F^\perp T^*TF^\perp\geqslant c^2 F^\perp$, we have $(E\wedge F^\perp)T^*T(E\wedge F^\perp)\geqslant c^2E\wedge F^\perp$.
  Since
  \begin{equation*}
      \|T(E\wedge F^\perp)\|\leqslant\|(T+A)(E\wedge F^\perp)\|+\|A(E\wedge F^\perp)\|\leqslant\|(T+A)E\|+\|A\|<c,
  \end{equation*}
  we get $E\wedge F^\perp=0$.
  This completes the proof.
\end{proof}

By \Cref{section 4}, we can derive the following conclusion with a proof analogous to that of \cite[Corollary 4.1]{Cui}.

\begin{corollary}\label{T-zI}
  Let $\mathcal{A}$ be an $AW^*$-algebra and $T\in\mathcal{K}_\mathcal{A}$.
  Then $zI+T$ is clean for any $z\in\mathbb{C}$.
  More precisely, there is an idempotent $P\in\mathcal{A}$ such that $zI+T-P$ is invertible and $\|(zI+T-P)^{-1}\|\leqslant 8$.
\end{corollary}

\begin{lemma}\label{exist V and E}
  Let $\mathcal{A}$ be a countably decomposable infinite $AW^*$-factor and $T\in\mathcal{A}$.
  Suppose $T-zI\notin\mathcal{K}_\mathcal{A}$ for all $z\in\mathbb{C}$ and the spectral projection of $|T|$ with respect to $[0,c]$ is finite for some $c\in(0,1)$.
  Then there is an invertible element $V$ and a projection $E$ in $\mathcal{A}$ satisfying the following conditions:
  \begin{enumerate}
    \item [$(1)$] $E\sim E^\perp$,
    \item [$(2)$] $\|EV^{-1}TVE\|<1$,
    \item [$(3)$] $E^\perp(V^{-1}TV)^*(V^{-1}TV)E^\perp\geqslant a^2 E^\perp$ for some $a>0$,
    \item [$(4)$] there is a finite projection $F\leqslant E$ such that $E(V^{-1}TV)E^\perp S=E-F$ for some $S\in E^\perp\mathcal{A}E$.
  \end{enumerate}
\end{lemma}

\begin{proof}
  By \Cref{T-cI compact}, there is a projection $E_0$ such that $E_0\sim E_0^\perp$ and $E_0^\perp T^*E_0\notin\mathcal{K}_\mathcal{A}$.
  By \Cref{compact equivalent}, the spectral projection $E_1$ of $|E_0^\perp T^*E_0|$ associated with $(a_0,+\infty)$ is infinite for some $a_0\in(0,1)$.
  It follows that $E_1\leqslant E_0$ and $E_1TE_0^\perp T^*E_1=E_1(E_0TE_0^\perp T^*E_0)E_1\geqslant a_0^2 E_1$.
  Then
  \begin{equation}\label{E1}
      E_1\sim E_1^\perp,\quad E_1TE_1^\perp T^*E_1\geqslant a_0^2 E_1.
  \end{equation}
  Let $A=E_1^\perp T^*S_0TE_1\in E_1^\perp\mathcal{A}E_1$, where $S_0$ is the inverse of $E_1TE_1^\perp T^*E_1$ in $E_1\mathcal{A}E_1$.
  Then we have $E_1TA=E_1TE_1$, i.e., $E_1T(I-A)E_1=0$.
  Let $b=\max\{1,2\|(I+A)T(I-A)\|\}$ and
  \begin{equation*}
      V=(I-A)(E_1+bE_1^\perp).
  \end{equation*}
  Since $A^2=0$, $V$ is invertible and $V^{-1}=(E_1+\frac{1}{b}E_1^\perp)(I+A)$.

  Let $T_0=V^{-1}TV$.
  Then $E_1T_0E_1=0$ and $\|T_0E_1\|=\|E_1^\perp T_0E_1\|\leqslant\frac{1}{2}$.
  By inequality \eqref{E1}, we have
  \begin{equation}\label{E1T0E1perp}
      E_1T_0E_1^\perp T_0^*E_1=b^2E_1TE_1^\perp T^*E_1\geqslant a_0^2b^2E_1.
  \end{equation}
  By \Cref{|T*|-|VT*V-1|}, the spectral projection $P$ of $|T_0|$ with respect to $[0,2a]$ is finite for some $a\in(0,\frac{1}{2})$.
  Let $E$ be the spectral projection of $|T_0E_1^\perp|$ with respect to $[0,a]$.
  Then we have $\|T_0(E-E_1)\|=\|T_0E_1^\perp E\|\leqslant a<2a$.
  By \Cref{TF<a}, we obtain $E_1^\perp-E^\perp=E-E_1\preceq P$.
  Since $P$ is finite and $E_1\sim E_1^\perp$, one can get $E\sim E^\perp$, i.e., condition $(1)$ holds.
  Moreover, we have
  \begin{equation*}
    \|ET_0E\|\leqslant\|T_0E_1\|+\|T_0(E-E_1)\|\leqslant\frac{1}{2}+a<1
  \end{equation*}
  and $E^\perp(E_1^\perp T_0^*T_0E_1^\perp)E^\perp\geqslant a^2E^\perp$, i.e., conditions (2) and (3) hold.

  By inequality \eqref{E1T0E1perp}, there exists $S_1$ such that $E_1T_0E_1^\perp S_1=E_1$.
  Since $E-E_1$ is a finite projection, we obtain $K=ET_0E^\perp S_1-E\in\mathcal{F}_\mathcal{\mathcal{A}}$.
  Let $F=R(KE)$ and $S=E^\perp S_1(E-F)$.
  It follows that $ET_0E^\perp S=(K+E)(E-F)=E-F$, i.e., condition $(4)$ holds.
\end{proof}

\begin{lemma}\label{infinite spectral proiection}
  Let $\mathcal{A}$ be a countably decomposable infinite $AW^*$-factor and $T\in\mathcal{A}$.
  If there is a projection $E$ satisfying the following two conditions:
  \begin{enumerate}
    \item [$(1)$] $E-ETE$ is invertible in $E\mathcal{A}E$ and $E^\perp T^*TE^\perp\geqslant c^2E^\perp$ for some $c>0$,
    \item [$(2)$] $L(EFE-E\wedge F)E_0$ is an infinite projection, where $F=I-L(TE^\perp)$ and $E_0$ is the spectral projection of $EFE-E\wedge F$ with respect to $[0,d]$ for some $d\in(0,1)$,
  \end{enumerate}
  then $T$ is clean.
  More precisely, there is an idempotent $P\in\mathcal{A}$ such that $T-P$ is invertible and
  \begin{equation}\label{(T-P')-1}
    \|(T-P)^{-1}\|\leqslant\left(1-\sqrt{\frac{8}{9-d}}\right)^{-\frac{1}{2}}\max\left\{\|(E-ETE)^{-1}\|, c^{-1}\right\}.
  \end{equation}
\end{lemma}

\begin{proof}
  We use the notation in \Cref{EFPQ}.
  Then $EFE-E\wedge F=E_{11}H$ and $L(EFE-E\wedge F)=E_{11}$.
  Let $P_0=E_{11}E_0$.
  Then $P_0$ is an infinite projection in $\mathcal{A}$ and we can write $P_0=P_1+P_2$, where $P_1$ and $P_2$ are infinite projections commuting with $E_{11}H$.
  Let $P_3=E_{11}-P_0$.
  For $k=1,2,3$, let
  \begin{equation*}
   I_k=E_{11}P_kE_{11}+E_{21}P_kE_{12},\quad Q_k=QI_k,\quad H_k=HI_k,\quad E_{ij}^{(k)}=E_{ij}I_k.
  \end{equation*}
  Let $E_1=E\wedge F^\perp$ and $E_2=E^\perp\wedge F$.
  Since $E_{22}^{(1)}$ is infinite, there exist partial isometries $V_1$ and $W_1'$ such that $V_1^*V_1=E_1$, $W_1'^*W_1'=E_{22}^{(1)}$ and $V_1V_1^*+W_1'W_1'^*=E_{22}^{(1)}$.
  Similarly, there exist partial isometries  $V_2$ and $W_2'$ in $\mathcal{A}$ such that $V_2V_2^*=E_2$, $W_2'W_2'^*=E_{11}^{(2)}$ and $V_2^*V_2+W_2'^*W_2'=E_{11}^{(2)}$.
  Let
  \begin{equation*}
      W_1=E_{12}^{(1)}W_1'E_{21}^{(1)}+E_{22}^{(1)}W_1'E_{22}^{(1)},\quad W_2=E_{11}^{(2)}W_2'E_{11}^{(2)}+E_{21}^{(2)}W_2'E_{12}^{(2)}.
  \end{equation*}
  It follows that $W_k$ commutes with $\{E_{ij}^{(k)}\}_{1\leqslant i,j\leqslant 2}$ for $k=1,2$, $W_1^*W_1=I_1$ and $W_2W_2^*=I_2$.
  Let
  \begin{equation*}
  \begin{split}
A&=a\left(V_1+E_{21}^{(1)}W_1+V_2+E_{21}^{(2)}W_2\right)+iE_{21}^{(3)}\in E^\perp\mathcal{A}E,\\
    T_1&=E_{11}^{(1)}+a\left(V_1+E_{21}^{(1)}W_1\right)+E_1,\\
    T_2&=E_{11}^{(2)}+a\left(V_2+E_{21}^{(2)}W_2\right),
  \end{split}
  \end{equation*}
  where $a=\frac{2\sqrt{2}}{\sqrt{1-d}}$.
  Then
  \begin{equation*}
      L(E+A)=L(T_1)+L(T_2)+E\wedge F+\frac{1}{2}\left(I_3+iE_{21}^{(3)}-iE_{12}^{(3)}\right).
  \end{equation*}
  By \Cref{idempotent}, there exists an idempotent $P\in\mathcal{A}$ such that $R(P)=E^\perp$, $P-E\in E^\perp\mathcal{A}E$ and $L((T-P)E)=L(E+A)$.

  To use \Cref{T invertible in A}, next we will prove that $L((T-P)E)-F^\perp$ is invertible in $\mathcal{A}$.
  We only need to prove $L(T_1)-(E_1+I_1-Q_1)$ is invertible in $(E_1+I_1)\mathcal{A}(E_1+I_1)$ and $L(T_2)-(I_2-Q_2)$ is invertible in $(E_2+I_2)\mathcal{A}(E_2+I_2)$.
  By the definition of $H_k$, we have $\|H_k\|\leqslant d<1$ for $k=1,2$.
  Therefore, $\sqrt{I_k-H_k}$ is invertible in $I_k\mathcal{A}I_k$ with inverse $\sqrt{I_k-H_k}^{-1}$.
  For $k=1,2$, let
  \begin{equation*}
      S_k=\sqrt{H_k}W_k^*+a\sqrt{I_k-H_k}=a\sqrt{I_k-H_k}\left(I_k+\frac{1}{a}\sqrt{I_k-H_k}^{-1}\sqrt{H_k}W_k^*\right).
  \end{equation*}
  Then $S_k$ commutes with $\{E_{ij}^{(k)}\}_{1\leqslant i,j\leqslant 2}$.
  Since $\|\frac{1}{a}\sqrt{I_k-H_k}^{-1}\sqrt{H_k}W_k^*\|<1$,
   $S_k$ is invertible in $I_k\mathcal{A}I_k$ with inverse $S_k^{-1}$.
  Let $U_k=E_{11}^{(k)}\sqrt{H_k}+E_{12}^{(k)}\sqrt{I_k-H_k}$.
  Then $U_k^*U_k=Q_k$ and $U_kU_k^*=E_{11}^{(k)}$.
  Since $W_1^*V_1=W_2V_2^*=0$, we obtain
  \begin{equation*}
      \begin{split}
           U_1 T_1&=E_{12}^{(1)}S_1\left(V_1+E_{21}^{(1)}W_1\right), \\
          (E_2+U_2)T_2&=(aE_2+S_2)\left(E_2+E_{11}^{(2)}+S_2^{-1}\sqrt{H_2}V_2^*\right)\left(V_2+E_{11}^{(2)}W_2\right).
      \end{split}
  \end{equation*}
  Let
  \begin{equation*}
  \begin{split}
      A_1&=\left(V_1+E_{21}^{(1)}W_1\right)^*S_1^{-1}E_{21}^{(1)}U_1,\\
      A_2&=\left(V_2+E_{11}^{(2)}W_2\right)^*\left(E_2+E_{11}^{(2)}-S_2^{-1}\sqrt{H_2}V_2^*\right)\left(\frac{1}{a}E_2+S_2^{-1}\right)(E_2+U_2).
  \end{split}
  \end{equation*}
  Then we have $Q_1 T_1A_1=Q_1$, $T_1A_1Q_1 T_1=T_1$, $(E_2+Q_2)T_2A_2=E_2+Q_2$ and $T_2A_2(E_2+Q_2)T_2=T_2$.
  Hence,
  \begin{equation}\label{F1LP(T')1}
  \begin{split}
      Q_1 L(T_1)T_1A_1&=Q_1,\quad T_1A_1Q_1 L(T_1)=L(T_1),\\
    (E_2+Q_2)L(T_2)T_2A_2&=E_2+Q_2,\quad T_2A_2(E_2+Q_2)L(T_2)=L(T_2).
  \end{split}
  \end{equation}
  Since
  \begin{equation*}
      \begin{split}
          Q_1&=L(Q_1 L(T_1))=Q_1-Q_1\wedge(E_1+I_1-L(T_1)),\\
          E_2+Q_2&=L((E_2+Q_2)L(T_2))=E_2+Q_2-(E_2+Q_2)\wedge(E_2+I_2-L(T_2)),
      \end{split}
  \end{equation*}
  we have $Q_1\wedge(E_1+I_1-L(T_1))=(E_2+Q_2)\wedge(E_2+I_2-L(T_2))=0$.
  Therefore, we can obtain $(E_1+I_1-Q_1)\vee L(T_1)=E_1+I_1$ and $(I_2-Q_2)\vee L(T_2)=E_2+I_2$.
  Combining with \Cref{E-F invertible} and \eqref{F1LP(T')1}, $L(T_1)-(E_1+I_1-Q_1)$ is invertible in $(E_1+I_1)\mathcal{A}(E_1+I_1)$ and $L(T_2)-(I_2-Q_2)$ is invertible in $(E_2+I_2)\mathcal{A}(E_2+I_2)$.

  With a similar proof of \Cref{Tclean}, $T-P$ satisfies $(a)$ and $(b)$ in \Cref{T invertible in A} for $a_1=c$ and $a_2=\frac{1}{\|(E-ETE)^{-1}\|}$.
  Therefore, $T-P$ is invertible.
  Next we consider the norm of $(T-P)^{-1}$.
  Since
  \begin{equation*}
  \begin{split}
       T_1T_1^*&=E_1+E_{11}^{(1)}+a^2 E_{22}^{(1)}+a\left(V_1+E_{21}^{(1)}W_1\right)+a\left(V_1+E_{21}^{(1)}W_1\right)^*,\\
       T_2T_2^*&=a^2E_2+E_{11}^{(2)}+a^2 E_{22}^{(2)}+a\left(V_2+E_{21}^{(2)}W_2\right)+a\left(V_2+E_{21}^{(2)}W_2\right)^*,
  \end{split}
  \end{equation*}
  we have $(T_kT_k^*)^2=(a^2+1)T_kT_k^*$ for $k=1,2$.
  Then $\|T_k\|=\sqrt{a^2+1}=\sqrt{\frac{9-d}{1-d}}$.
  We get $\|A_1\|\leqslant\|S_1^{-1}\|\leqslant 1$ and $\|A_2\|\leqslant\left(1+\sqrt{d}\|S_2^{-1}\|\right)\sqrt{\frac{1}{a^2}+\|S_2^{-1}\|^2}\leqslant 1$ since $\|S_k^{-1}\|\leqslant\frac{1}{2\sqrt{2}-\sqrt{d}}$.
  Hence, $\|T_kA_k\|\leqslant\sqrt{\frac{9-d}{1-d}}$.
  It follows from \eqref{F1LP(T')1} that $L(T_1)\leqslant \|T_1A_1\|^2 L(T_1)Q_1 L(T_1)$ and $L(T_2)\leqslant\|T_2A_2\|^2 L(T_2)(E_2+Q_2)L(T_2)$.
  Then $L(T_1)(E_1+I_1-Q_1)L(T_1)\leqslant\frac{8}{9-d}L(T_1)$, $L(T_2)(I_2-Q_2)L(T_2)\leqslant\frac{8}{9-d}L(T_2)$.
  Therefore, $\|(E_1+I_1-Q_1)L(T_1)\|\leqslant\sqrt{\frac{8}{9-d}}$ and $\|(I_2-Q_2)L(T_2)\|\leqslant\sqrt{\frac{8}{9-d}}$.
  By \Cref{P vee Q}, $\|L((T-P')E)L((T-P')E^\perp)\|\leqslant \sqrt{\frac{8}{9-d}}$.
  Then we obtain inequality \eqref{(T-P')-1} by \Cref{T invertible in A}.
  This completes the proof.
\end{proof}

\begin{lemma}\label{compact and spectral projection equivalent}
  Let $E$ be a projection in an infinite $AW^*$-factor $\mathcal{A}$ and $T\in\mathcal{A}$.
  Suppose $F=I-L(TE^\perp)$.
  Then $EF^\perp E-E\wedge F^\perp\notin\mathcal{K}_\mathcal{A}$ if and only if $L(EFE-E\wedge F)E_0$ is an infinite projection, where $E_0$ is the spectral projection of $EFE-E\wedge F$ with respect to $[0,d]$ for some $d\in(0,1)$.
\end{lemma}

\begin{proof}
  With notation in \Cref{EFPQ}, then $EFE-E\wedge F=E_{11}H$ and $EF^\perp E-E\wedge F^\perp=E_{11}-E_{11}H$.
  Hence $L(EFE-E\wedge F)=E_{11}$.
  Since $E_{11}E_0$ is equal to the spectral projection of $E_{11}-E_{11}H$ with respect to $[1-d,1]$, we complete the proof by \Cref{compact equivalent}.
\end{proof}

Let $T$ and $S$ be elements in an $AW^*$-algebra $\mathcal{A}$.
We say that $T$ and $S$ are \emph{similar} if there is an invertible element $V\in\mathcal{A}$ such that $S=V^{-1}TV$.
By \Cref{Tclean}, \Cref{infinite spectral proiection} and \Cref{compact and spectral projection equivalent}, we can obtain the following conclusion.

\begin{corollary}\label{infinite clean}
  Let $T$ be an element in a countably decomposable infinite $AW^*$-factor $\mathcal{A}$.
  If there is a projection $E\in\mathcal{A}$ such that
  \begin{itemize}
    \item[$(1)$] $E-ETE$ is invertible in $E\mathcal{A}E$ and $E^\perp T^*TE^\perp\geqslant c^2 E^\perp$ for some $c>0$,
    \item[$(2)$] $E\wedge L(TE^\perp)\sim E^\perp\wedge(I-L(TE^\perp))$ or $EL(TE^\perp)E-E\wedge L(TE^\perp)\notin\mathcal{K}_\mathcal{A}$,
  \end{itemize}
  then all elements similar to $T$ are clean.
\end{corollary}

\begin{lemma}\label{P and Pperp infinite}
Let $T$ be an element in a countably decomposable infinite $AW^*$-factor $\mathcal{A}$.
If the spectral projections of $|T|$ with respect to $[0,c]$ and $(c,+\infty)$ are both infinite for every $c\in(0,1)$, then for every $a>0$, there exists an invertible element $V$ and a projection $E$  satisfying the following conditions:
   \begin{enumerate}
     \item [$(1)$] $E\sim E^\perp$,
     \item [$(2)$] $\|V-I\|\leqslant a$,
     \item [$(3)$] $\|EV^{-1}TVE\|\leqslant\frac{3}{4}$,
    \item [$(4)$] $E^\perp(V^{-1}TV)^*V^{-1}TVE^\perp\geqslant\frac{1}{16}E^\perp$,
     \item [$(5)$] $EV^{-1}TVE^\perp$ is not compact relative to $\mathcal{A}$,
     \item [$(6)$] For any $\varepsilon>0$, there exists $E_0\leqslant E$ such that $E_0\sim E-E_0$ and $\|V^{-1}TVE_0\|\leqslant\varepsilon$.
   \end{enumerate}
\end{lemma}

\begin{proof}
  Let $E$ be the spectral projection of $|T|$ with respect to $\big[0,\frac{1}{2}\big]$.
  If $ETE^\perp\notin\mathcal{K}_\mathcal{A}$, then let $V=I$ and we complete the proof.
  Next we assume that $ETE^\perp\in\mathcal{K}_\mathcal{A}$.
  By \Cref{compact equivalent}, there is a finite projection $E_1\leqslant E^\perp$ such that $\|ET(E^\perp-E_1)\|\leqslant\frac{1}{4}$.
  Since $E^\perp T^*TE^\perp\geqslant\frac{1}{4}E^\perp$, we have
  \begin{equation}\label{(I-E)T}
      (E^\perp-E_1)T^*E^\perp T(E^\perp-E_1)\geqslant\frac{1}{16}(E^\perp-E_1).
  \end{equation}
  Let $E_2$ be the spectral projection of $|T|$ with respect to $\big[0,\frac{1}{8}\big]$.
  Then $E_2$ is an infinite subprojection of $E$ and $\|TE_2\|\leqslant\frac{1}{8}$.
  Let $F=W^*E_2 W\wedge (E^\perp-E_1)$, where $W$ is a partial isometry such that $W^*W=E^\perp$ and $WW^*=E$.
  Then $F$ is an infinite projection in $\mathcal{A}$.
  Since $\|TE_2\|\leqslant\frac{1}{8}$, by inequality \eqref{(I-E)T} and \Cref{A-B}, $F(WTE^\perp-ETW)^*(WTE^\perp-ETW)F\geqslant\frac{1}{64}F$.
  Therefore, $WTE^\perp-ETW\notin\mathcal{K}_\mathcal{A}$.
  It follows that $(WTE^\perp-ETW)+tWTW\notin\mathcal{K}_\mathcal{A}$ for some $0<t<\min\left\{a,\frac{1}{16\|T\|}\right\}$.
  Note that $(I+tW)^{-1}=I-tW$ and we have
  \begin{equation*}
      E(I-tW)T(I+tW)E^\perp=ETE^\perp-t(WTE^\perp-ETW)-t^2WTW\notin\mathcal{K}_{\mathcal{A}}.
  \end{equation*}
  Let $V=I+tW$.
  Then $V$ and $E$ satisfy conditions $(1)-(6)$.
\end{proof}

\begin{lemma}\label{EFE-E wedge F compact}
  Let $E$ be a projection in an infinite $AW^*$-factor $\mathcal{A}$ and $T\in\mathcal{A}$.
  If $F=L(TE^\perp)$ and $EFE-E\wedge F\in\mathcal{K}_\mathcal{A}$, then $(E-E\wedge F)TE^\perp\in\mathcal{K}_\mathcal{A}$.
\end{lemma}

\begin{proof}
  We use the notation in \Cref{EFPQ}.
  Then
  \begin{align*}
   (E-E\wedge F)TE^\perp&=(E-E\wedge F)FTE^\perp=\left(E_{11}H+E_{12}\sqrt{H(I_0-H)}\right)TE^\perp\\
   &=E_{11}\sqrt{H}\left(\sqrt{H}+E_{12}\sqrt{(I_0-H)}\right)TE^\perp.
  \end{align*}
  Since $E_{11}\sqrt{H}=\sqrt{EFE-E\wedge F}\in\mathcal{K}_{\mathcal{A}}$, we have $(E-E\wedge F)TE^\perp\in\mathcal{K}_\mathcal{A}$.
\end{proof}

\begin{lemma}\label{contradiction}
  Let $T$ be an element and $E$ a projection in a countably decomposable infinite $AW^*$-factor $\mathcal{A}$.
  Let $F=L(TE^\perp)$.
  Suppose the following conditions are satisfied:
  \begin{enumerate}
    \item [$(a)$] $E\sim E^\perp$,
    \item [$(b)$] $E^\perp T^*TE^\perp\geqslant c^2E^\perp$ for some $c>0$,
    \item [$(c)$] $E\wedge F\not\sim E^\perp\wedge F^\perp$,
    \item [$(d)$] $EFE-E\wedge F$ is compact relative to $\mathcal{A}$.
  \end{enumerate}
  Then we have the following conclusions:
  \begin{enumerate}
    \item [$(1)$] If $ETE^\perp\notin\mathcal{K}_\mathcal{A}$, then $E\wedge F$ is infinite and there is a finite projection $E_1\leqslant E^\perp$ such that $E^\perp-E_1=E^\perp TE^\perp S_1$ for some $S_1\in\mathcal{A}$.
    \item [$(2)$] If there exists a finite projection $E_2\leqslant E$ such that $E-E_2=ETE^\perp S_2$ for some $S_2\in\mathcal{A}$, then the spectral projection of $|T^*|$ with respect to $[0,c_1]$ is finite for some $c_1\in(0,1)$ and the spectral projection of $|T|$ with respect to $[0,\varepsilon]$ is infinite for any $\varepsilon>0$.
  \end{enumerate}
\end{lemma}

\begin{proof}
  By condition $(d)$ and \Cref{EFE-E wedge F compact}, $(E-E\wedge F)TE^\perp\in\mathcal{K}_\mathcal{A}$.

  (1) Since $ETE^\perp\notin\mathcal{K}_\mathcal{A}$ and $(E-E\wedge F)TE^\perp\in\mathcal{K}_\mathcal{A}$, we can obtain $(E\wedge F)TE^\perp\notin\mathcal{K}_\mathcal{A}$.
  Therefore, $E\wedge F$ is infinite.
  By condition $(b)$ and \Cref{TES}, there exists $S\in\mathcal{A}$ such that $TE^\perp S=L(TE^\perp)=F$.
  Then $TE^\perp S(F-E\wedge F)=F-E\wedge F$.
  Let $P=L(E^\perp S(F-E\wedge F))\leqslant E^\perp$.
  Then
  \begin{equation}\label{L(TP)}
    L(TP)=L(TE^\perp S(F-E\wedge F))=F-E\wedge F.
  \end{equation}
  Since $(E\wedge F)TP=0$, we obtain $ETP=(E-E\wedge F)TE^\perp P\in\mathcal{K}_{\mathcal{A}}$.
  By \Cref{compact equivalent}, we have $(P-P_1)T^*ET(P-P_1)\leqslant \frac{c^2}{2}(P-P_1)$ for some finite projection $P_1\leqslant P$.
  Then by condition $(b)$, we get
  \begin{equation}\label{P-P1}
    (P-P_1)T^*E^\perp T(P-P_1)\geqslant\frac{c^2}{2}(P-P_1).
  \end{equation}
  Since $E\wedge F$ is infinite, $E^\perp\wedge F^\perp$ is finite by condition $(c)$.
  By equation \eqref{L(TP)},
  \begin{equation*}
    E^\perp-L(E^\perp TP)=E^\perp-L(E^\perp L(TP))=E^\perp\wedge L(TP)^\perp=E^\perp\wedge F^\perp.
  \end{equation*}
  Then $E^\perp-L(E^\perp TP)$ is finite.
  Since $P_1$ is finite and
  \begin{align*}
    & \ \ \ \ L(E^\perp TP)-L(E^\perp T(P-P_1))\\
    &=L(E^\perp TP_1)\vee L(E^\perp T(P-P_1))-L(E^\perp T(P-P_1))&\\
          &\sim L(E^\perp TP_1)-L(E^\perp TP_1)\wedge L(E^\perp T(P-P_1)\\
          &\leqslant L(E^\perp TP_1)\sim R(E^\perp TP_1)\leqslant P_1,
  \end{align*}
  the projection $L(E^\perp TP)-L(E^\perp T(P-P_1))$ is finite.
  Let
  \begin{equation*}
    E_1=E^\perp-L(E^\perp T(P-P_1))=E^\perp-L(E^\perp TP)+L(E^\perp TP)-L(E^\perp T(P-P_1)).
  \end{equation*}
  Then $E_1$ is a finite subprojection of $E^\perp$.
  By \Cref{TES} and inequality \eqref{P-P1}, there exists $S_0$ such that $E^\perp T(P-P_1)S_0=L(E^\perp T(P-P_1))=E^\perp-E_1$.
  Let $S_1=(P-P_1)S_0$.
  Then $E^\perp TE^\perp S_1=E^\perp-E_1$.

  (2) Let $W$ be a partial isometry such that $W^*W=E$ and $WW^*=E^\perp$.
  Then $(E-E\wedge F)(ETE^\perp+W)=(E-E\wedge F)TE^\perp\in\mathcal{K}_\mathcal{A}$.
  Moreover,
  \begin{equation*}
    (ETE^\perp+W)(E^\perp S_2+W^*)=ETE^\perp S_2+WW^*=E-E_2+E^\perp=E_2^\perp.
  \end{equation*}
  By \Cref{ET compact}, $E-E\wedge F$ is finite.
  Since $E$ is infinite, $E\wedge F$ is infinite.
  Hence $E^\perp\wedge F^\perp$ is finite by condition $(c)$.
  Since
  \begin{equation*}
    E-E\wedge F\sim E\vee F-F=F^\perp-E^\perp\wedge F^\perp,
  \end{equation*}
  $F^\perp$ is finite.
  By condition $(b)$ and \Cref{TES}, there exists $S\in\mathcal{A}$ such that $TE^\perp S=L(TE^\perp)=F$.
  Then by \Cref{|T*|-|VT*V-1|}, the spectral projection of $|T^*|$ with respect to $[0,c_1]$ is finite for some $c_1\in(0,1)$.

  For any $\varepsilon>0$, let $P_\varepsilon$ be the spectral projection of $|T|$ with respect to $[0,\varepsilon]$.
  Then we have $P_\varepsilon^\perp T^*TP_\varepsilon^\perp\geqslant \varepsilon^2 P_\varepsilon^\perp$.
  We claim that there exists an infinite projection $Q\in\mathcal{A}$ such that $\|TQ\|<\varepsilon$.
  In this case, we have $Q\wedge P_\varepsilon^\perp=0$.
  Then $Q\preceq P_\varepsilon$ and hence $P_\varepsilon$ is infinite.
  Next we prove the claim.
  If $TE\in\mathcal{K}_\mathcal{A}$, then we can find such $Q$ by \Cref{compact equivalent}.
  If $TE\notin\mathcal{K}_\mathcal{A}$, then by \Cref{compact equivalent}, the spectral projection $P$ of $|TE|$ with respect to $(c_2,+\infty)$ is infinite for some $c_2\in(0,1)$.
  We have $P\leqslant E$ and $PT^*TP\geqslant c_2^2 P$.
  By \Cref{TES}, there exists $S_3\in\mathcal{A}$ such that $TPS_3=L(TP)$.
  Recall that $TE^\perp S=F$.
  Let $A=E^\perp S-PS_3$ and $B=A(F\wedge L(TP))$.
  Then $TE^\perp B=F\wedge L(TP)$.
  Hence $F\wedge L(TP)\leqslant R(B)$.
  Since $F^\perp$ is finite and $L(TP)$ is infinite, $F\wedge L(TP)$ is infinite.
  Then $R(B)$ is infinite and hence $L(B)$ is infinite.
  Since $TB=0$, we have $TL(B)=0$.
  Then $L(B)\leqslant R(T)^\perp$ and hence $R(T)^\perp$ is infinite.
  Let $Q=R(T)^\perp$.
  Then $TQ=0$.
  We complete the proof.
\end{proof}


Next we prove our main result in this section.

\begin{theorem}\label{section 5}
    All countably decomposable infinite $AW^*$-factors are clean.
\end{theorem}

\begin{proof}
  Since $T$ is clean if and only if $T^*$ is clean, we only need to consider the following two cases:
  \begin{enumerate}
    \item [\textbf{I.}] There exists $c\in(0,1)$ such that the spectral projection of $|T|$ with respect to $[0,c]$ or $(c,+\infty)$ is finite.
    \item [\textbf{II.}] For any $c\in(0,1)$, the spectral projections of $|T|$ and $|T^*|$ with respect to $[0,c]$ and $(c,+\infty)$ are infinite.
  \end{enumerate}

  \textbf{Case I.} Since $T$ is clean if and only if $I-T$ is clean, by \Cref{T and I-T}, we can assume that the spectral projection of $|T|$ with respect to $[0,c]$ is finite.
  If there exists $z\in\mathbb{C}$ such that $T-zI\in\mathcal{K}_\mathcal{A}$, then $T$ is clean by \Cref{T-zI}.
  If  $T-zI\notin\mathcal{K}_\mathcal{A}$ for any $z\in\mathbb{C}$, then there is an invertible element $V$ and a projection $E$ in $\mathcal{A}$ satisfying conditions (1)-(4) of \Cref{exist V and E}.
  By \Cref{|T*|-|VT*V-1|}, the spectral projection of $|V^{-1}TV|$ with respect to $[0,c_V]$ is finite for some $c_V\in(0,1)$.
  It follows from \Cref{contradiction}(2) that the condition $(c)$ or $(d)$ of \Cref{contradiction} is not satisfied for $V^{-1}TV$.
  Then $V^{-1}TV$ and $T$ are both clean by \Cref{infinite clean}.

  \textbf{Case II.} By \Cref{P and Pperp infinite}, we can assume there exists a projection $E\in\mathcal{A}$ such that $E\sim E^\perp$ and
  \begin{enumerate}
     \item [$(a)$] $\|ETE\|\leqslant\frac{3}{4}$,
     \item [$(b)$] $E^\perp T^*TE^\perp\geqslant\frac{1}{16}E^\perp$.
     \item [$(c)$] $ETE^\perp\notin\mathcal{K}_\mathcal{A}$,
     \item [$(d)$] for any $\varepsilon>0$, there exists $E_0\leqslant E$ with $E_0\sim E-E_0$ and $\|TE_0\|\leqslant\varepsilon$.
  \end{enumerate}
  Let $F=L(TE^\perp)$.
  If $E\wedge F\sim E^\perp\wedge F^\perp$ or $EFE-E\wedge F\notin\mathcal{K}_\mathcal{A}$, then $T$ is clean by \Cref{infinite clean}.
  Next we assume $E\wedge F\not\sim E^\perp\wedge F^\perp$ and $EFE-E\wedge F\in\mathcal{K}_\mathcal{A}$.

  It follows \Cref{contradiction}(1) that $E\wedge F$ is infinite and there exists a finite projection $E_1\leqslant E^\perp$ such that $E^\perp-E_1=E^\perp TE^\perp S$ for some $S\in\mathcal{A}$.
  Then $E^\perp\wedge F^\perp$ is finite.
  Similar to the proof of equation \eqref{L(TP)}, there exists a projection $P\leqslant E^\perp$ such that $L(TP)=F-E\wedge F\sim F\vee E-E=E^\perp-E^\perp\wedge F^\perp$.
  Therefore, $L(TP)$ is infinite and hence $P$ is infinite.
  Since
  \begin{align*}
    E\wedge F=F-L(TP)&=L(T(E^\perp-P))\vee L(TP)-L(TP)\\
    &\sim L(T(E^\perp-P))-L(T(E^\perp-P))\wedge L(TP),
  \end{align*}
  the projection $L(T(E^\perp-P))$ is infinite.
  Therefore, $E^\perp-P$ is infinite and $P\sim E^\perp-P$.
  By condition $(d)$, there exists $E_0\leqslant E$ such that $E_0\sim E-E_0$ and $\|TE_0\|\leqslant\frac{c_0}{4}$.
  Let $W=W_1+W_2$, where $W_1$ and $W_2$ are partial isometries in $\mathcal{A}$ such that
  \begin{equation*}
    W_1^*W_1=P,\quad W_1W_1^*=E_0,\quad  W_2^*W_2=E^\perp-P,\quad W_2 W_2^*=E-E_0.
  \end{equation*}
  Then $W^*W=E^\perp$ and $WW^*=E$.
  Let $V=I+tW$ with $t\in\left(0,\frac{1}{1+24\|T\|}\right)$.
  Then $V^{-1}=I-tW$.
  Let $T_1=V^{-1}TV=T-t(WT-TW+tWTW)$ and $B=T-T_1$.
  Then $\|BE^\perp\|\leqslant\|B\|\leqslant\frac{1}{8}$.
  Since $E^\perp T^*TE^\perp\geqslant\frac{1}{16}E^\perp$, by \Cref{A-B}, we have $E^\perp T_1^*T_1 E^\perp=E^\perp(T-B)^*(T-B)E^\perp\geqslant\frac{1}{64}E^\perp$.
  Since $\|ET_1E\|\leqslant \|ETE\|+t\|WTE\|\leqslant\frac{7}{8}$, $E-ET_1E$ is invertible in $E\mathcal{A}E$.
  Since $E^\perp T^*TE^\perp\geqslant\frac{1}{16}E^\perp$ and $P\leqslant E^\perp$, by \Cref{TES}, we have $TPS_1=L(TP)$ for some $S_1\in\mathcal{A}$.
  Then
  \begin{equation*}
    E^\perp TPS_1TE^\perp=E^\perp L(TP)TE^\perp=E^\perp FTE^\perp=E^\perp TE^\perp.
  \end{equation*}
  Let $S_2=S_1 TE^\perp S$.
  Then $E^\perp TPS_2=E^\perp TE^\perp S=E^\perp-E_1$.
  By \Cref{|T*|-|VT*V-1|}, the spectral projection of $|(E^\perp TP)^*|$ in the subalgebra $E^\perp\mathcal{A}E^\perp$ with respect to $[0,c_0]$ is finite for some  $c_0\in(0,1)$.
  Then the spectral projection of $|(tWTP)^*|$ in $E\mathcal{A}E$ with respect to $[0,tc_0]$ is finite.
  Let $T_2=ETP$, $T_3=EBP=t(WTP-ETWP+tWTWP)$ and $A=t(tWTWP-ETWP)$.
  Then $\|A\|\leqslant 2t\|TWP\|=2t\|TE_0W_1\|\leqslant\frac{tc_0}{2}$.
  Let $P_1$ be the spectral projection of $|T_3^*|$ in $E\mathcal{A}E$ with respest to $[0,\frac{tc_0}{4}]$.
  Then we have
  \begin{equation*}
    \|((tWTP)^*+A^*)P_1\|=\|T_3^*P_1\|\leqslant\frac{tc_0}{4}<tc_0-\|A\|.
  \end{equation*}
  By \Cref{TF<a}, $P_1$ is finite.

  Since $EFE-E\wedge F\in\mathcal{K}_\mathcal{A}$, by \Cref{EFE-E wedge F compact}, $T_2=(E-E\wedge F)TP\in\mathcal{K}_\mathcal{A}$.
  By \Cref{compact equivalent}, there exists a projection $Q\in P\mathcal{A}P$ such that $P-Q$ is finite and $\|T_2Q\|\leqslant\frac{tc_0}{8}$.
  Let $P_2$ be the spectral projection of $|(T_2Q-T_3)^*|$ in $E\mathcal{A}E$ with respest to $[0,\frac{tc_0}{16}]$.
  Then $\|(T_3^*-(T_2Q)^*)P_2\|\leqslant\frac{tc_0}{16}<\frac{tc_0}{4}-\|T_2Q\|$.
  By \Cref{TF<a}, we have $P_2\preceq P_1$ and hence $P_2$ is finite.
  Then by \Cref{|T*|-|VT*V-1|}, there exists a finite projection $E_2\leqslant E$ such that $(T_2Q-T_3)S_3=E-E_2$ for some $S_3\in\mathcal{A}$.
  Since $P-Q$ is finite,  we have $T_2(P-Q)S_3\in\mathcal{F}_\mathcal{A}$.
  Therefore, $E_3=E_2\vee R(T_2(P-Q)S_3)$ is finite and
  \begin{equation*}
  \begin{split}
      ET_1PS_3(E-E_3)&=(T_2-T_3)S_3(E-E_3)\\
      &=(T_2Q-T_3)S_3(E-E_3)+T_2(P-Q)S_3(E-E_3)=E-E_3.
  \end{split}
  \end{equation*}
  We have $E\wedge L(T_1E^\perp)\sim E^\perp\wedge L(T_1E^\perp)^\perp$ or $EL(T_1E^\perp)E-E\wedge L(T_1E^\perp)\notin\mathcal{K}_\mathcal{A}$ by \Cref{|T*|-|VT*V-1|} and \Cref{contradiction}(2).
  Then $T_1$ and $T$ are both clean by \Cref{infinite clean}.
\end{proof}

We end this section with the following problem.

\begin{question}
    Are all infinite $AW^*$-factors clean?
    More generally, are all $AW^*$-algebras clean?
\end{question}



\end{document}